\documentclass[%
 reprint,
onecolumn,
 amsmath,amssymb, aps,
]{revtex4-2}

\usepackage{dcolumn}
\usepackage{bm}
\usepackage{color}
\usepackage{amstext}
\usepackage{amssymb}
\usepackage{graphicx}
\usepackage{amsmath}
\usepackage{tcolorbox}
\usepackage{nccmath}
\usepackage{amsfonts}
\usepackage[utf8]{inputenc}
\usepackage{comment}
\usepackage{siunitx}
\numberwithin{equation}{section}

\DeclareMathAlphabet\mathbfcal{OMS}{cmsy}{b}{n}

\newcommand\const{\mathrm{const}}

\newcommand\vv{\boldsymbol{v}}

\begin{document}

\title{{\Large Mikhail Alekseevich Lavrentiev and Fluid Dynamics.\\[4mm] On the 125th anniversary of his birth }}

\author{ V.A. Vladimirov}
\affiliation{Universities of York and Leeds, UK}

\author{Translated by O.M. Lavrentieva}
\affiliation{Technion - Israel Institute of Technology}

%\date{January 2024}
\maketitle
\Large
\vskip -2mm
\textit{Abstract:}
The purpose of the present paper is to advertise and illuminate the applied research of \textit{Mikhail Alekseyevich Lavrentiev (1900–1980)}, a distinguished scientist who generated new ideas and obtained many outstanding results in mechanics and mathematics.
He also initiated and maintained numerous flourishing activities in mechanics and applied mathematics in the former Soviet Union.
 He was notably an outstanding researcher in pure and applied mathematics as well as in engineering.
On the applied side, his research spanned fluid dynamics, explosion physics, mechanics, technology, and several branches of engineering.
It was always aimed at the practical challenges of his time and achieved impressive successes.  It has been fascinating to uncover all of this.  Much of his biographical material and a review of his activities in pure mathematics have been published, but no review of his applied achievements is available.  This paper aims to give Lavrentiev's applied research the recognition it deserves in the communities of fluid dynamics, applied mathematics, engineering, and the history of science.

\textit{Keywords:} Lavrentiev, fluid mechanics, hydrodynamics, explosion physics, engineering, technology, wing theory, jet theory, separated flows, hydrodynamic behaviour of materials, cumulation, lined charges, explosive welding, wave theory, solitons, vortex dynamics, locomotion of snake and fish.

\textit{Caption of the Figure}.
The paper's Contents and meaning frame are presented below as a map scheme, showing M.A. Lavrentiev's applied research results.
The central block shows the unity between the cumulation (or focusing) of an explosion's impact and the hydrodynamic behaviour of materials.
The big (pink) circles identify the main areas of research in which he has personally obtained results.
Smaller (green) circles represent secondary research tasks, sometimes initiated only by Lavrentiev and developed with his support and participation.
The principles that connect different topics are shown in the (yellow) smoothed rectangles.
The arrows propose the logical and historical interconnections between different topics in his research.
Some terms coined by Lavrentiev (\textit{e.g.} “pestle”) are clarified within the text.
\begin{center}
      \resizebox{17cm}{25cm}{\includegraphics[angle=90]{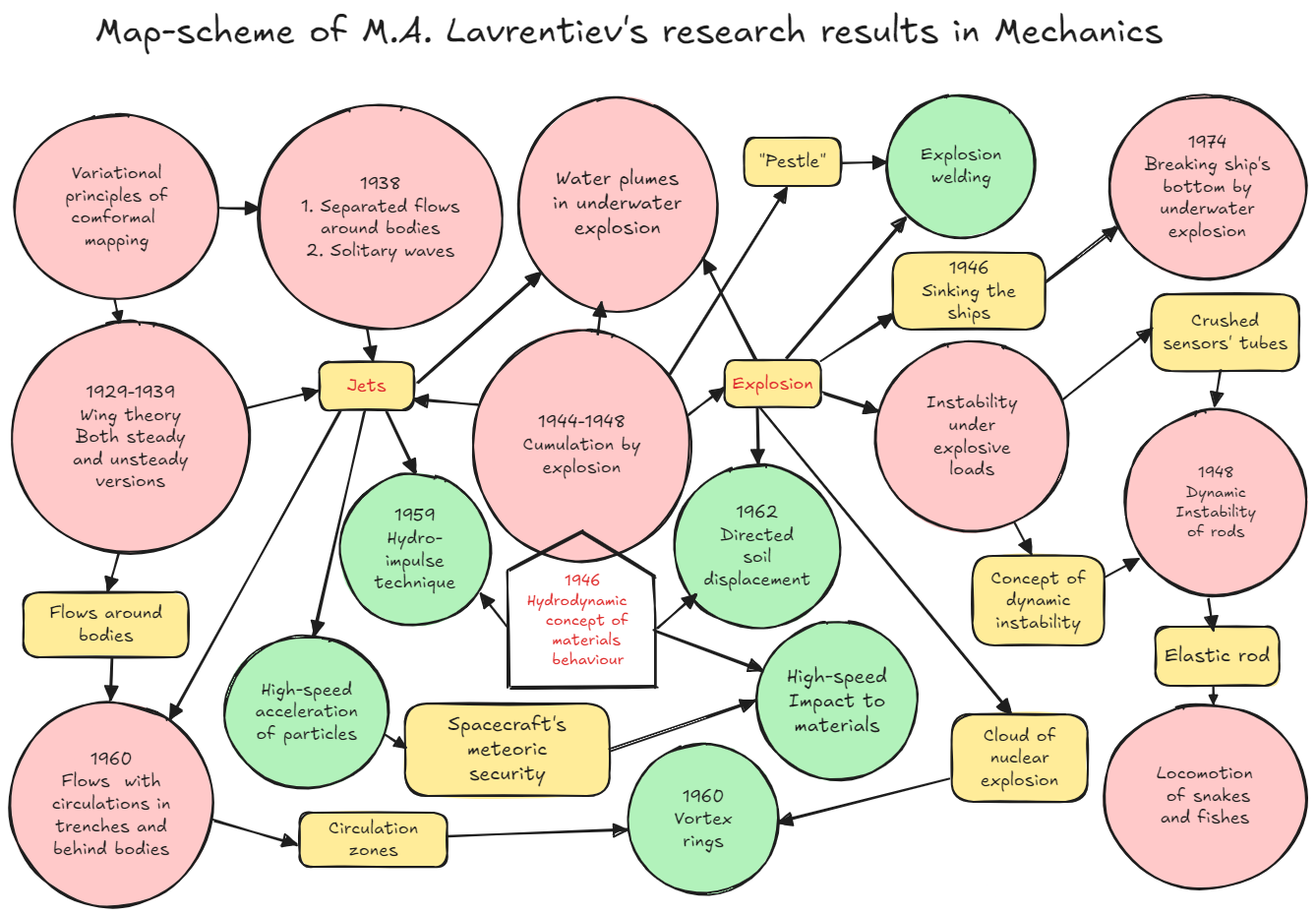}}
      %{etest.pdf}}
    %\captionof{figure}{Your caption here}
    % \label{fig:your-label}
\end{center}

\section{Introduction}
"Hydrodynamics is one of the most ancient sciences.  Today (as often before), it is experiencing \textit{one more period of its youth}."
These words open Lavrentiev's book \emph {“Problems of Hydrodynamics and their Mathematical Models”} [1].
He explained the rapid development of hydrodynamics in the 1960s and 1970s through the widespread adoption of computers, the essential extension of mathematical methods, and the rapid increase of practical problems requiring urgent attention.
The expression “one more period of its youth” was deliberately used.
In the 1930s, Lavrentiev had participated in an earlier period of intense hydrodynamic research activity driven by advances in aviation. From the beginning of the 20th century, many relevant studies in fluid mechanics were focused on the problem of flow around a solid body, particularly aerofoils.
The Moscow school of hydrodynamics, headed by S.A. Chaplygin, made significant contributions to the wing theory.
Their achievements included the development of models for flow over the wing and the extension of the limits of applicability of these models.
They built on the fundamental results of N.E. Joukowski, L. Prandtl, M.V. Kutta, and S.A. Chaplygin, who had determined the mechanism of lift and had found exact analytical solutions for special wing profiles.
The theoretical group in the Central Aero-Hydrodynamic Institute (TsAGI), located near Moscow, included M.V. Keldysh, A.P. Chaplygin, and L.V. Sedov, as well as S.A. Chaplygin and Lavrentiev.
Their investigations involved precise problem formulation, advanced mathematical methods, and clear recognition that the objective was application to engineering practice.
Lavrentiev was already well-known as a pure mathematician when he was engaged to TsAGI in 1929.
The following text describes his \emph {applied} results and their various interconnections.
Finally, we summarise his views on the objectives of applied research.
Several historical remarks at the end of the paper were added during the translation.

\section{Applied research Directions and results}

\textit{1. Theory of steady wings}.
\vspace{1mm}

 From 1929 to 1934,  Lavrentiev solved several important problems in wing theory and rapidly became known as an active researcher in aerodynamics and fluid mechanics.
 In a paper from this period [2], Lavrentiev developed a general method for solving the flow problem around a thin profile (an arc) of an arbitrary shape,
using mathematical techniques that were new for these problems.
These included reducing the problem to a singular integral equation of the first kind, developing an algorithm to construct solutions, and proving its convergence.
After extending the wing theory to thin profiles of arbitrary shape, he determined the best profile in this class.
Using the variational properties of conformal mappings, he showed that a circular arc maximises the lift force [3].
\vspace{3mm}

\textit{2. Theory of unsteady wings and hydrodynamic impact.}
\vspace{1mm}

During the same period, Lavrentiev completed three important tasks in collaboration with his student and later well-known Soviet academician, M.V. Keldysh.
The first was devoted to the theory of oscillating wings [4].
The authors considered the flow of an incompressible, inviscid fluid over a thin, slightly curved profile that performs prescribed oscillations.
For the first time, the theory of functions of complex variables was used to describe unsteady fluid motion.
This led to the determination of the time-dependent components of forces and torques, thereby generalising the classical results of Joukowski and Chaplygin.
In particular, oscillation regimes with additional traction have been described [4].
Later, these results were applied to the theory of aerodynamic flutter.
In [5], the motion of a wing profile under the free surface of a liquid was studied for the first time, using
 thin-wing and small-amplitude wave approximations as the basis of the theory.
The corresponding linearised problem was reduced to an integral equation of the first kind.
The authors noted that the singular integral operator in this equation transforms polynomials into polynomials and used this remarkable fact to construct an approximate solution.
As a result, several simple analytic formulae for the forces exerted on the wing and, importantly, the wave resistance were obtained.
In the third joint paper [6], the authors obtained firm results in the theory of body impact on water,
where they found the general solution to a two-dimensional impact problem, which
 considered a body with curvilinear boundaries.
Lavrentiev became interested in this problem while modelling a seaplane water landing.
\vspace{3mm}

\textit{3. Separated Flows and Solitary Waves.}
\vspace{1mm}

Lavrentiev succeeded in developing the wing theory through his elegant use of integral equations and complex function theory.
At the same time, aerodynamic problems stimulated him to develop conformal mapping theory and, in particular, the variational principles of conformal mappings.
The latter fell within his range while solving the problem of optimisation of the wing shape.
He proposed modifications to these variational principles and used them to prove several existence and uniqueness theorems for various flows of inviscid incompressible fluid over bodies with flow separation.
While searching for further applications for these methods, Lavrentiev also turned his attention to wave problems.
 Here, in 1943-1947, he finally proved the existence of a solution to the Euler equations (of inviscid fluid dynamics) describing the propagation of a solitary surface wave [7,8].
 Such a wave had been famously observed in 1844 by J.~Scott Russell.
Boussinesq and Rayleigh independently obtained approximate descriptions of this phenomenon.
 Lavrentiev proved the existence of a family of steady periodic waves of arbitrarily long period $T$ using his conformal mapping and variational methods [7,8].
 In the limit $T\to\infty$,  a symmetric wave with a single maximum was obtained.
Within the framework of his approach, the Boussinesq and Rayleigh results turned out to be equivalent to  approximating the surface velocity $\vv$ by the Lavrentiev equation
\begin{equation}
    |\vv|^2=(1+2f f''/3)/f^2,\nonumber
\end{equation}
where $y=f(x)$ is the equation of the liquid surface.
By combining this formula with the Bernoulli integral $|\vv|^2+2gf=\const.$,  Lavrentiev derived the simplest mathematical model for a solitary wave retaining only quadratic terms, and
the Boussinesq and Rayleigh theories were thereby validated.
Lavrentiev's work promoted further developments in non-linear wave theory and related mathematical methods.
\vspace{3mm}

\textit{4. Cumulation (concentrated impact) of the explosion.}
\vspace{1mm}

The Second World War was a crucial period in Lavrentiev's applied research.
With the clarity characteristic of an outstanding mathematician, he turned to military and engineering problems.
His investigations in 1944-1948 led to the creation of the theory of cumulation by explosion [9-12].
This phenomenon was discovered in the 1880s.
It consisted of localising the impact (including momentum and energy) exerted by an explosive charge on a barrier when the charge cavity had an opening on the side facing the barrier.
This local effect was significantly increased if the explosive cavity was lined with a thin metal layer.
In this case, the explosion penetrated the barrier to considerable depth.
World War II spurred the exploitation of this phenomenon in the design of new weapons.
However, engineering projects in this field were based solely on empirical studies and were carried out independently (in highly classified research) in different countries.
Many researchers, including G.I.Pokrovsky in the USSR, tried to construct a theory of this effect.
The theory proposed by Lavrentiev was based on two new ideas [1,12,13].
The main assumption is that any solid material (and particularly a metal) under intense stress behaves essentially like a liquid.
Lavrentiev further understood that the internal strength forces of a material can be neglected in comparison with the pressure differences and stresses induced by an explosion or high-speed collision.
Thus, the corresponding motion of a solid material can, in a first approximation, be described as the flow of an \emph {inviscid} liquid.
This approach seems natural and obvious today, but considering a metal as a liquid was controversial at the time and sparked a heated debate.
The hydrodynamic analysis of the cumulation phenomenon led Lavrentiev to the second important idea:
that an active element breaking a barrier can be modelled by a high-speed, cumulative jet formed from the cavity lining.
As a result, the construction of the theory of cumulation was split into modelling two processes: jet formation and its interaction with a barrier.
Using image analysis, both processes were formulated as variants of the same problem, the collision of two liquid jets.
Two versions of the theory were constructed: a plane model and an axisymmetric model.
The collision problem of two plane jets was solved first, and Lavrentiev used this to build the theory.
For axisymmetric jets, Lavrentiev successfully used conservation laws (for the fluxes of mass, momentum, and energy) to derive the main relations governing the process.
The most important prediction of his theory was his formula for the depth $L$ at which the jet penetrates a barrier:
 \begin{equation}
    L=l\sqrt{\rho_1/\rho_2},\nonumber
\end{equation}
where $l$ denotes the length of the jet between the source and impact, and $\rho_1$ and $\rho_2$ are the densities of the jet and barrier materials, respectively.
Between 1944 and 1948, Lavrentiev conducted experiments on the cumulation theory in Feofanija, near Kyiv.
These confirmed the validity of this formula, revealing that $L$ is independent of the speed of the jet.
It is worth mentioning that G. Birkhoff, G.I. Taylor, and their collaborators independently developed a similar theory of cumulation by explosion.
Later, Lavrentiev cited their 1948 paper as the first open (unclassified) publication on this problem.
This period marked a highly productive phase in Lavrentiev’s career.
Alongside his primary accomplishments, he also obtained several results, whose true significance became apparent later.
\vspace{3mm}

\textit{5. Heuristic concept of cumulation and its multiple outcomes.}
\vspace{1mm}

Lavrentiev introduced the general heuristic concept of cumulation in a continuous medium as any process that increases energy density in small volumes at the expense of a decrease in large ones, and analysed many phenomena from this point of view.
After 1940, the connection between two subjects, fluid mechanics and the physics of explosions,  became the central area of his applied research.
The hydrodynamic concept of material behaviour allowed him not only to explain the phenomenon of cumulation by explosion but also to use his powerful intuition to model various practical problems.
This led to such applications as the hydro-impulse technique, the acceleration of solid particles to meteoric speeds, and explosive welding.
Other examples include explosion-induced collapse of cylindrical and spherical shells and free surface jet flows (plumes) generated by underwater explosions.
Less well-known are the cumulative effects during the propagation of a tsunami along an underwater mountain ridge [1] and the focusing of a surface wave induced by an underwater explosion of a ring or semi-ring-shaped charge [1,12,13].
Surprisingly, he even suggested using cumulation to refloat ships that had sunk at sea.
\vspace{3mm}

\textit{6. Directional displacement of soil by explosion.}
\vspace{1mm}

The solution of the problem of explosion-induced directional soil displacement was a significant theoretical and practical success of the hydrodynamic treatment of material behaviour under impulse loading.
This problem  had long existed;
in particular, directional explosions were used during dam construction.
However, the directivity was incomplete and, after the explosion, a considerable portion of the soil was blown away from the main direction.
Lavrentiev formulated the objective in its simplest form:
\textit{the selected soil volume should fly towards a required destination as a solid body.}
The basic elements of his model are: a) the soil behaves as an ideal incompressible fluid;
b) an explosive layer continuously covers the selected soil volume, and
c) the momentum transmitted to the soil (per unit of area) is proportional to the thickness of that layer.
An elegant mathematical solution to this problem, in which the velocity potential is a linear function of the coordinates, is given in [9].
This solution requires that the thickness of the explosive layer decrease linearly in the direction of throw.
Special experimental explosions confirmed the validity of the solution and later enabled corrections to be made to the compressibility and mechanical strength of the medium.
These corrections did not change the formulation of the problem, but they enabled a more accurate determination of the positions of the charges and the order of their ignition to achieve the planned objective.
An impressive practical implementation of these ideas was the construction of the gigantic Medeu Dam, which protects the city of Almaty, Kazakhstan.
The dam was constructed by a series of five directed blasts, which displaced five million cubic metres of rock.
The massive final explosion (3600 tonnes of an ammonium nitrate-based explosive) occurred on 21 October 1966.
It was used to collapse the walls of the gorge and form the bulk of the barrier.
A subsequent reinforcement explosion occurred in 1967, and the dam reached its functional design profile in 1972.
The engineering proved its worth almost immediately: on July 15, 1973, it successfully held back a catastrophic mudflow that would otherwise have repeated the devastation of 1921, which killed around 500 people and destroyed a large part of the city.
Lavrentiev played a pivotal role in the project as the primary scientist and leading consultant behind the ``directed blast" method
(see also https://www.youtube.com/watch?v=DmeqVCN12JA).

\vspace{3mm}

\textit{7. Hydro-impulse guns.}
\vspace{1mm}

The hydro-impulse technique is another important application that stems from Lavrentiev's work on cumulation.
To verify the hydrodynamic theory of jet penetration in a material, he proposed testing it with a water jet impinging on a layer of coal, which serves as the medium to be broken.
In 1948, an experiment was constructed that became a prototype of modern hydro-cannons.
In 1959-66, Lavrentiev returned to the problem of obtaining a high-speed water jet capable of breaking firm materials, particularly coal layers.
Under his leadership, B.V.Vojcekhovski at the Institute of Hydrodynamics (Novosibirsk) constructed several hydro-impulse set-ups.
The speeds achieved by the water jets were close to 1 km/s and, with the use of the so-called ``cavitators", exceeded 3 km/s [13].
A jet of this kind penetrates metals to a considerable depth (e.g. in copper, up to 100 mm).
\vspace{3mm}

\textit{8. Protecting spacecraft from meteor impact}
\vspace{1mm}

Exploitation of the cumulation approach led to progress in several other areas of technical physics.
With the high-speed cumulative jet  at their disposal, Lavrentiev and his disciple V.M. Titov used it to accelerate solid particles to high speeds [13].
At the Institute of Hydrodynamics, they developed a method to accelerate particles using a jet of detonation products from a cylindrical cumulative explosive charge.
The millimetre-size particles were accelerated up to 12-15 km/s.
Such speeds are required to study the impact of meteors on spacecraft.
\vspace{3mm}

\textit{9. Explosive welding.}
\vspace{1mm}

From 1944-48,  Lavrentiev conducted multiple experiments on the cumulation-of-explosion theory in Feofanija near Kyiv.
This was a period of extremely fruitful work; see [13] and the introduction to [12].
In addition to his main achievements, several important secondary results were obtained during this period.
For example, the first series of experiments included charges with a double-layer cavity lining, where the remarkable phenomenon of explosion welding was observed.
After an explosion, a high-speed cumulative jet and a low-speed massive body (Lavrentiev called it a ``pestle”) consisting of both lining metals were formed.
The pestle was a monolithic body in which the initially separated metals were welded over the entire contact surface.
Concurrently, Lavrentiev's co-worker N.M.Sytyj obtained a monolithic copper bar by compressing a pack of copper wires by an explosion.
However, the importance of these results and similar publications in the USA was not appreciated at that time.
A surge of interest in this new method of producing multilayer materials arose in the 1960s, when intensive studies of explosion welding, stimulated by technological applications, began in several laboratories worldwide.
This research and technology direction was initiated by Lavrentiev and carried out by his disciple, A.A.Deribas and colleagues.
\vspace{3mm}

\textit{10. Programme to sink ships.}
\vspace{1mm}

In 1946, M.A. Lavrentiev participated in a government programme to sink trophy German ships.
This work initiated his interest in underwater explosions.
In 1947, he launched a series of experiments in Feofanija to rupture a ship's hull by underwater explosions and obtained a paradoxical result.
The purpose of the experiments was to find the minimal charge, $q$, needed for the rupture as a function of its distance $h$ from the hull.
Intuitively, the function $q(h)$ should increase monotonically.
In practice,  $q$ was constant over a certain range of $h$, i.e., hull rupture due to more distant charges did not require an increase in the weight of charge.
Lavrentiev analysed the effect and concluded that it is caused by a cumulative jet directed to the hull.
The jet forms when a cavity caused by the detonation product collapses.
Short descriptions of these results can be found in [1,13].
\vspace{3mm}

\textit{11. Explosion plumes on the free surface of water.}
\vspace{1mm}

Much later, Lavrentiev studied one more effect related to explosion jets.
In 1969, he described problems concerning the formation and structure of plumes on a free surface observed after an underwater explosion.
Many experimental investigations of plumes were carried out in the 1940s.
However, no satisfactory theory of the phenomenon had been given.
Lavrentiev suggested a model in which the key element was a depression on the free surface that acts as a cumulation cavity.
This depression is formed by the shock wave from the explosion, which cuts off and spreads part of the liquid above the charge, leaving a depression on the free surface.
 The surface cavity then collapses due to the fluid motion induced by the rapid extension of the underwater cavity filled with detonation products.
The cumulative jet then forms a plume.
 Computation of the water motion in 1971 showed that this model predicts an upward jet.
However, a broad set of experiments in the field and laboratory, and extensive computational work led by Lavrentiev, showed that the mechanism of plume formation varies with the depth of the charge location.
It appears that the plume always represents a cumulative jet, but the main mechanism of this jet formation is not always connected with a depression on the free surface.
These studies present a few cases of these mechanisms and provide detailed descriptions of each [1,13].
\vspace{3mm}

\textit{12. Dynamic instability.}
\vspace{1mm}

During experiments on explosions in Feofanija, Lavrentiev observed an unusual deformation of the submerged pipes that mounted the sensors.
Further research and analysis of this effect led to the theory of dynamic instability of elastic bodies.
The simplest example of this effect is the bending instability of an elastic rod under an impulsively applied longitudinal load, which was examined in a joint paper with the  A.Yu.Ishlinsky [10], later a known Soviet academician.
Lavrentiev played a pivotal role in Ishlinsky's career as a senior colleague and patron.
The most impressive result of [10] was the finding of a principal difference in the behaviour of the rod under static and dynamic loads.
In the first case, only the rod shape with the longest wavelength is realised; in the second case, the result can differ.
Theoretically and experimentally, it was shown that if the longitudinal load exceeds $n$ times the critical Euler force, then the rod bending is sinusoidal with the number of semi-waves $m$ equal to the integer closest to $\sqrt{n/2}$.
If the rod is deformed in an elastic-plastic way without destruction, its final shape is close to sinusoidal with the same semi-wave number $m$.
If it breaks up, the number of fractures is also $m$.
This work combines an understanding of the phenomenon with the simplicity of the result, and professionals appreciate it as a classic.
\vspace{3mm}

\textit{13. Locomotion of snakes and fish.}
\vspace{1mm}

The key elements of Lavrentiev's scientific considerations became part of his research arsenal and generated new, often unexpected, analogues.
People around him often thought that he needed much less information than they did to obtain meaningful results.
He preferred to formulate a hypothesis about the main mechanism of a phenomenon, temporarily ignoring the stages of its empirical investigation and the step-by-step accumulation of knowledge.
Due to his keen intuition and focus on analysing the qualitative aspects of phenomena, his hypotheses often proved true, and surprised or stimulated other researchers.
For example, the paper [11], coauthored with his son, M.M.Lavrentiev, offers an unexpected analogy.
It is devoted to explaining the self-propelled motion of living organisms such as fish and a grass snake.
As a model, the authors considered the motion of a flexible elastic rod inside a smooth, rigid channel.
Note that the key element of the model here is, once again, an elastic rod.
The energy source for the motion is the muscle force distribution of an animal that bends its body.
The total force exerted on the rod comes from the counter-pressure of the channel's walls.
The model may describe the motion of a snake in the grass or  swimming in water.
Hydrodynamic pressure or stem reactions serve as channel walls for motion in water and grass, respectively.
This unusual model has interested many researchers, and discussions on the problem continue.
\vspace{3mm}

\textit{14. Vortex flows and vortex rings.}
\vspace{1mm}

From the early 1960s, Lavrentiev paid great attention to the study of vortex flows.
His interest was sparked by the government's requests to study two problems: the feasibility of storing waste in ocean-bottom depressions and to describe the dynamics of a nuclear explosion cloud.
This initiated his interest in vortex rings and he was fascinated by their unusual properties.
His public and teaching lectures often included an account of them, which he accompanied with an experimental demonstration of smoke rings.
Lavrentiev's favourite method of formulating problems was to lead them to a paradox or, at least, to a contrast with another well-known phenomenon.
 He contrasted the extended propagation of a vortex ring with the rapid stopping of a child's air balloon of the same volume.
He saw this difference as an opportunity to considerably reduce hydrodynamic drag.
While developing general ideas concerning vortex rings and discussing them with his disciples, Lavrentiev raised the problems of describing their formation mechanisms, their structure for a given formation process, the law of motion, the maximum distance travelled by a ring, and the amount of admixture that a vortex can transfer.
The experimental programmes at the Institute of Hydrodynamics in Novosibirsk addressed most of these problems in detail.
These experiments were carried out over a wide range of Reynolds numbers for vortex rings (from $10^2$ to $10^7$) in air and water.
Based on the results, Lavrentiev's disciple B.A.Lugovtsov developed a self-similar theory of turbulent vortex rings, using vortex momentum as the sole characteristic parameter to define the universal structure of vortex rings at very high Reynolds numbers, thereby partially completing the research stage [1,13].

\vspace{3mm}
\textit{15. Trench flow modes. Flow separation and circulation zones.}
\vspace{1mm}

The government project on waste disposal in ocean depressions (e.g., ocean trenches) prompted  Lavrentiev to construct mathematical models of flow around bodies with flow separation and a rear circulation zone.
In 1962, he considered such flows in plane geometry, modelling the circulation zone as a point vortex or a domain with constant vorticity [1,13]; the latter model
 (independently of G.K.Batchelor) was accepted as the main working option.
Important problems considered by Lavrentiev’s collaborators included flows around a cylinder or plate, flows over a trench or step at the bottom, and flows along a straight wall with an adjoint closed circulation zone (a plane analogue of Hill's vortex) [1].
The experimental programmes initiated by Lavrentiev indicated that waste disposal will be rapidly transferred to the external flow above the depression, and so did not support the original proposal.
However, these studies led Lavrentiev to return to his earlier studies of flow around bodies [1].
His interest in the subject dates back to his work at TsAGI and his 1938 publication on existence theorems, and he maintained this interest throughout his life.

\section{ Miscellaneous and general}

\textit{Naturalist's approach.}
\vspace{1mm}

Being a mathematician, Lavrentiev also often demonstrated his ability as a naturalist.
As examples, one may consider the following suggestions that were widely discussed in the USSR.
He faced the problem of understanding the mechanism of swimming fish in 1954 when he used a ring-shaped pool in Katsiveli,
Crimea [13].
The experiments included studies of dolphin swimming modes.
Lavrentiev actively participated in the observations, and they initiated or accelerated an entire spectrum of his further ideas and studies.
In addition  to  explosion welding and fish and grass-snake motions, these included his renewed interests in dynamic instability, underwater explosions, and water waves.
One can also mention his surprising hypotheses on the damping of waves in reservoirs by rain and the nature of the impulsive strong wind called the “Novorossiysk Bora,” where, according to his suggestion, the air circulation zone formed behind the near-coast mounting crest could spontaneously jump above it and damage the coastal city of Novorossiysk on the Black Sea.
\vspace{3mm}

\textit{Setting up research projects and programmes.}
\vspace{1mm}

After 1957, Lavrentiev's research activity changed significantly.
He was crucially involved in the setting up and development of research programmes for the newly established  Siberian Division of the USSR Academy of Sciences.
This means that he worked with a huge amount of information and was occupied with formulating or criticising hypotheses concerning various phenomena.
However, young and bright researchers around him were inexperienced in the most challenging tasks of choosing and formulating research problems.
He formulated interesting and important problems for many of his students, disciples, and collaborators.
At the same time, for mature researchers, he was the greatest authority to assess their prospects and validate the results obtained.
This activity peaked in 1960-1961, when the Institute of Hydrodynamics, the first of the Siberian Division Institutes, started up under his leadership.
\vspace{3mm}

\textit{The key role of the monograph [1].}
\vspace{2mm}

The monograph [1] summarises Lavrentiev's activity in mechanics and considers the problems he was directly or indirectly involved in.
Many of these problems are mentioned above, and the list indicates his key role in shaping both the conceptions and the subject matter.
The monograph has no analogues in the literature.
Perhaps it can be compared with the well-known book “Hydrodynamics” by Garrett Birkhoff.
It includes paradoxes, concentrates on the qualitative side of the phenomena, and tends to provide the simplest descriptions of the material.
The emphasis on the qualitative side of phenomena and the simplicity of their description is not occasional but rather reflects Lavrentiev's understanding of the aim of research in mechanics.
The understanding of what should be seen as a result of mechanics was so important for him that he addressed it in the preface to [1], which we cite here in length due to its importance for the understanding of this paper:

“In the previous century and at the beginning of ours (he meant the 19th and 20th centuries, VV),  books in hydrodynamics consisted of long calculations using elementary and special functions.
The modern American researcher S. Goldstein wrote that it was impossible to imagine that the water considered in these calculations was wet.
Today, many papers also present lengthy, complex results based on exact solutions of differential equations that are far removed from reality.
The practical value of these works is relatively low because the hydrodynamic equations themselves are quite approximate when some important physical phenomena are concerned.
That’s why some results in the so-called accurate theory involve calculations with many digits for values that represent only rough approximations of the exact ones.
This book contains no complicated calculations or sophisticated theories.
The most interesting physical processes are so complicated that in the contemporary state of science, constructing a universal theory that applies across all stages of the process is rarely possible.
Instead, we need to reveal, through experiments and observations, the leading factors that control certain stages of the process.
Having selected these factors, we ignore the less important factors and construct the simplest mathematical scheme (a model of the process), considering only the selected factors.
Sometimes, corrections should be made that consider minor but essential factors.
This can be done using additional algorithms that can be applied to solve model problems.
To obtain the general scheme of the process, one needs to combine solutions of several local problems.
This is performed with the help of some general ideas, such as velocity field continuity, and so on.”

Lavrentiev consistently applied this understanding of the objective of mechanics to practice.
He always insisted that scientific results be obtained using the simplest methods and the fewest resources.
He extended this principle to both theory and experiment.
The work on accelerating solid particles with a cumulative gas jet [13] illustrates this principle.
From the beginning of this project, Lavrentiev encouraged his students to find simple, fast, and achievable ways of acceleration instead of developing traditional methods (the so-called light gas guns).
\vspace{3mm}

\textit{How was it all possible?}
\vspace{1mm}

One may ask: how could a mathematician be so practical in understanding the objectives of mechanics and consistently implement these objectives?
Perhaps the reason was Lavrentiev's unique breadth and dynamism, and his opposition to all artificial barriers and limitations.
It is characteristic that his main achievement, the theory of cumulation by explosion, came from breaking the barriers between hydrodynamics and the physics of explosion.
Undoubtedly,  Lavrentiev's achievements have a long-lasting influence in several branches of mechanics..
Many of his results are included in the professional and popular literature.
Several of his radical hypotheses have generated discussion, and thus stimulated scientific progress to date.
\vskip 3mm

Acknowledgements: The author is grateful to Dr O.M.Lavrentieva (Technion - Israel Institute of Technology), the granddaughter of M.A.Lavrentiev, for her hard, devoted, and enthusiastic work on this paper. The authors thank Professor H.K. Moffatt for editing this paper.

\vskip 6mm

\textbf{REFERENCES:}
\vskip 1mm
$\circ$ Russian references are translated into English; for the originals see [15].
\vskip 1mm
$\circ$ TsAGI =  Central Aerohydrodynamic Institute, Moscow region, USSR.
\vskip 2mm

1. Lavrentiev, M.A. and Shabat, B.V. 1977 \textit{Problems of hydrodynamics and their mathematical models}, 2nd edition, Moscow: Nauka, 416p; (in Russian).
French translation: Lavrentiev, M. and Shabat, B. 1980 \textit{Effets Hydrodynamiques et modeles mathématiques.} Moscow: Mir, 365 p. https://ia601600.us.archive.org/31/items/\\m.-lavrentiev-
b.-chabat-effets-hydrodynamiques-et-modeles-mathematiques-mir\\-1980/M.

2. Lavrentiev, M.A. 1932 On constructing a flow around a given shape arc, \textit{Proceedings of TsAGI}, issue 118, 53p.; (in Russian).

3. Lavrentiev, M.A. 1934 On an extremum problem in aircraft wing theory, \textit{Proceedings of TsAGI}, issue 155, 40p.; (in Russian).

4. Keldysh, M.V. and Lavrentiev, M.A. 1935 On the theory of oscillating wings, \textit{TsAGI Technical Notes}, 45, 48-52; (in Russian).

5. Lavrentiev, M.A. and  Keldysh, M.V. 1937 On the motion of a wing beneath the surface of a heavy liquid,
 \textit{Proceedings of the Conference on Wave Resistance Theory}, Moscow, 31-64; (in Russian).

6. Lavrentiev, M.A. and Keldysh, M.V. 1935 The general problem of the hydrodynamic impact of a rigid body to water, \textit{Proceedings of CAHI}, issue 152, 5–12; (in Russian).

7. Lavrentiev, M.A. 1943 On the theory of long waves, \textit{Dokl. Acad. Nauk SSSR,} \textbf{41}, 275-277; (in Russian).
English translation: \textit{Am. Math. Soc. Transl.}, 1954, \textbf{102}, 51–53.

8. Lavrentiev, M.A. 1947 Towards the theories of long waves, \textit{Proceeding of Mathematical Institute, Ukraine}, 8, 13-69 (in Ukrainian);
see also: Lavrent'ev, M. A. 1954 Part I. On the theory of long waves; Part II. A contribution to the theory of long waves,  \textit{American Mathematical Society Translations}, \textbf{102}, Providence, 53p.

9.	Kuznetsov, V.M., Lavrentiev, M.A., and Sher, E.N. 1960 On directed soil displacement using explosive substances,
Journal of Applied Mechanics and Technical Physics, 1960, 4, 49-50; (in Russian).

10.	 Lavrentiev, M.A. and Ishlinsky, A.Yu. 1949 Dynamic forms of losing the stability of elastic systems, \textit{Doklady Akademii Nauk SSSR (Proceedings of the USSR Academy of Sciences),} \textbf{64}, 6, 779-782; (in Russian).

11.	Lavrentiev, M.A. and Lavrentiev, M.M. 1962 On the principle of creating a propulsive force for motion, \textit{Journal of Applied Mechanics and Technical Physics}, 4, 3-9; (in Russian); see also  Lavrentiev, M.A. 1973 Mathematical modeling of fishes and grass snakes motions. \textit{Journal of Applied Mechanics and Technical Physics}, 2, 164-165; (in Russian).

12. \textit{Selected Works by M.A. Lavrentiev: Mathematics and Mechanics}, 1990, Nauka, Moscow, 600 p; (in Russian).

13.  Lavrentiev, M.A. 2000 Experiences of Life. 50 Years in Science. In: \textit{The Age of Lavrentiev, Dobretsov, N. L. and Ibragimova, Z. M. (eds.)}, Novosibirsk: Siberian division of the Russian Academy of Science Publishing House, 456 p; (in Russian).
www.prometeus.nsc.ru/elibrary/2000vek/

14. Vladimirov, V.A. 1990 Mikhail Alekseevich Lavrentiev and hydrodynamics (on the 90th anniversary of his birth), \textit{Journ. Appl. Mech. Tech. Physics}, 6, 3-12, (in Russian).

15. Vakulenko, L.D. et al. 1985 \textit{Lavrentiev, M.A. 
(1900-1980): \\ Bibliographic index (ed.  Soboleva, E.B.)}, Novosibirsk: Institute of Hydrodynamics.
\\The internet address: https://www.prometeus.nsc.ru/akademgorodok/ \\ lavrentev/biblio/works1.ssi or\\
https://web.archive.org/web/20070807005021\\
/http://www.prometeus.nsc.ru/akademgorodok/lavrentev/biblio/works1.ssi

16. Papadopoulos, A. 2020 A note about Mikhaïl Lavrentieff and his world of analysis in the Soviet Union. In: \textit{Handbook of Teichmüller Theory, Volume VII (ed. A. Papadopoulos)}, Berlin: European Mathematical Society Publishing House, 317-347.

\vspace{10mm}

\section{APPENDIX, added after translating the paper:}
\vspace{3mm}

\textit{The status of this paper.}
\vspace{1mm}

This paper [14], published only in Russian in 1990, had an official status.
The Institute of Hydrodynamics and the Journal of Applied Mechanics and Technical Physics commissioned it in 1988 for the 90th anniversary of M.A. Lavrentiev's birth.
The text of the paper was also used in the Introduction to \textit{Selected Works by M.A. Lavrentiev} [12]; see the acknowledgements by the editor, Academician L.V.Ovsjannikov.
During the writing of this paper, the author and translator were affiliated with the Institute of Hydrodynamics in Novosibirsk, now renamed the Lavrentiev Institute of Hydrodynamics.
This paper is devoted to the 125th birthday of a great scientist, Mikhail Alekseevich Lavrentiev, whose legacy belongs to the whole World.
An additional reason to publish this paper is that the original Russian version appears not to be available, having been excluded from the official journal's Russian website without the author's permission.

\vspace{3mm}

\textit{Adaptations to the text in the English version.}
\vspace{1mm}

The abstract and the division into topical paragraphs with subtitles were added to modernise the paper's composition.
The additional references [12-16] were intended to compensate for the limited number of references in [14].
The book [13] contains Memoirs by M.A. Lavrentiev and memoir articles by his disciples.
The best bibliographic index of his publications is [15].
A review of his pure mathematical research in English is available in [16].
\vspace{3mm}

\textit{Influence of barriers between West and East.}
\vspace{1mm}

The lack of knowledge about Lavrentiev's applied results in the West is worth considering.
This may be due to several factors.
During Stalin's rule and the Cold War (say, 1930-1985), communications between the USSR and Western academic communities were severely suppressed.
The USSR dominated the Communist bloc of Eastern countries, in which the "Capitalist Western World" was considered the primary ideological, political, and military enemy.
Communism was considered the future of all mankind, and the USSR's sciences and technologies were often praised as superior.
Consequently, Russian was the dominant language in the Eastern academic communities, and publishing papers in the West was often considered "servility or humiliation before the West".
Lavrentiev belonged to the top of the USSR academic community and was the leading organiser of science; therefore, he needed to follow those restrictions and rules to continue those important roles.
Probably, it was the main factor explaining the absence of his original applied publications in English.
There were several officially published translations from internal USSR journals; the readership of such translations was very limited.
Notice that before 1930, he published  12 papers (all in French) in pure mathematics in central European journals; see the references in [15].
The second key restrictive factor was the practical nature of his applied results, often involving patents or classified information.
An additional factor (especially for internet users) is the existence of multiple transliterations of his name, including Lavrent’ev, Lavrentev, Lavrentieff, Lavrentjev, Lavrentyev, and Lavrentiev.

 \end{document}